\input lanlmac
\def\href#1#2{{#2}}

\input epsf.tex

\overfullrule=0mm

\newcount\figno
\figno=0
\def\fig#1#2#3{
\par\begingroup\parindent=0pt\leftskip=1cm\rightskip=1cm\parindent=0pt
\baselineskip=11pt
\global\advance\figno by 1
\midinsert
\epsfxsize=#3
\centerline{\epsfbox{#2}}
\vskip 12pt
{\bf Fig.\ \the\figno:} #1\par
\endinsert\endgroup\par
}
\def\figlabel#1{\xdef#1{\the\figno}}
\def\encadremath#1{\vbox{\hrule\hbox{\vrule\kern8pt\vbox{\kern8pt
\hbox{$\displaystyle #1$}\kern8pt}
\kern8pt\vrule}\hrule}}


\magnification=\magstep1
\baselineskip=12pt
\hsize=6.3truein
\vsize=8.7truein
 at 8truept
 at 8truept
 at 10truept

\vbox{\hfill IPhT-t13/076}

\bigskip\bigskip

\font\bigrm=cmr12 at 14pt \centerline{\bigrm A note on irreducible maps with several boundaries} 
 
\bigskip\bigskip

\centerline{J. Bouttier$^{1,2}$ and E. Guitter$^1$}
  \smallskip
\centerline{$^1$ Institut de Physique Th\'eorique}
  \centerline{CEA, IPhT, F-91191 Gif-sur-Yvette, France}
  \centerline{CNRS, URA 2306}
\centerline{$^2$ D\'epartement de Math\'ematiques et Applications}
\centerline{\'Ecole normale sup\'erieure, 45 rue d'Ulm, F-75231 Paris Cedex 05}  
\centerline{\tt jeremie.bouttier@cea.fr}
\centerline{\tt emmanuel.guitter@cea.fr}

  \bigskip

     \bigskip\bigskip

     \centerline{\bf Abstract}
     \smallskip
     {\narrower\noindent

We derive a formula for the generating function of $d$-irreducible bipartite planar maps with
several boundaries, i.e.\ having several marked faces of controlled degrees.
It extends a formula due to Collet and Fusy for the case of arbitrary (non necessarily irreducible) bipartite planar maps, which we recover by taking $d=0$.
As an application, we obtain an expression for the number of $d$-irreducible bipartite planar maps with a prescribed 
number of faces of each allowed degree. Very explicit expressions are given in the case of maps without multiple 
edges ($d=2$), $4$-irreducible maps and maps of girth at least $6$ ($d=4$). Our derivation is based 
on a tree interpretation of the various encountered generating functions.
\par}

     \bigskip

\nref\IRRED{J. Bouttier and E. Guitter, {\it On irreducible maps and slices}, arXiv:1303.3728 [math.CO].}
\nref\CF{G. Collet and \'E. Fusy, {\it A simple formula for the series of bipartite and quasi-bipartite maps 
with boundaries}, FPSAC 2012, Nagoya, Japan, DMTCS proc. {\bf AR} (2012) 607-618, arXiv:1205.5215 [math.CO].}
\nref\BFa{O. Bernardi and \'E. Fusy, {\it A bijection for triangulations, quadrangulations, pentagulations, etc.},
J.\ Combin.\ Theory Ser.\ A {\bf 119} (2012) 218Ð244, arXiv:1007.1292 [math.CO].}
\nref\BFb{O. Bernardi and \'E. Fusy, {\it Unified bijections for maps with prescribed degrees and girth},
J.\ Combin.\ Theory Ser.\ A {\bf 119} (2012) 1351Ð1387, arXiv:1102.3619 [math.CO].}
\nref\TutteCPM{W.T. Tutte, {\it A Census of planar maps}, Canad.\ J.\ Math.\ 
{\bf 15} (1963) 249-271.}
\nref\BrownNS{W.G. Brown, {\it Enumeration of non-separable planar maps}, Can.\ J.\ Math.\ {\bf 15} (1963) 526-554.}
\nref\BT{W.G. Brown and W.T. Tutte, {\it On the enumeration of rooted non-separable planar maps}, 
Canad.\ J.\ Math.\ {\bf 16} (1964), 572-577.}
\nref\SJ{G. Schaeffer and B. Jacquard, {\it A bijective census of non-separable planar maps}, J.\ Comb.\ Theory {\bf A 83} (1998) 1Ð20.}
\nref\DSQ{J. Bouttier and E. Guitter, {\it Distance statistics in quadrangulations with no multiple edges 
and the geometry of minbus}, J.\ Phys.\ A: Math.\ Theor.\ {\bf 43} (2010) 205207,
arXiv:1002.2552 [math-ph].}
\nref\OEIS{see Sequence A179300 in The On-Line Encyclopedia of Integer Sequences, published electronically at http://oeis.org, 2010.}
\nref\Schae{G. Schaeffer, {\it Conjugaison d'arbres
et cartes combinatoires al\'eatoires}, PhD Thesis, Universit\'e 
Bordeaux I (1998).}
\nref\SchaeBCEM{G. Schaeffer, {\it Bijective census and random generation of Eulerian planar maps}, Electron.\ J.\ Combin.\ {\bf 4} (1997)
R20.}
\nref\TutteCS{W.T. Tutte, {\it A Census of slicings},
Canad.\ J.\ Math. {\bf 14} (1962) 708-722.}

\newsec{Introduction}

This note is an extension of a preceding paper \IRRED\ on the enumeration of 
$d$-irreducible planar maps via slice decomposition. We shall use here a number of results from this
paper, which the reader is invited to consult for explicit proofs.

We recall that a {\it planar map} is a cellular embedding of a graph in the sphere, considered up to continuous deformation.
A map is {\it bipartite} if one can color its vertices in black and white so that no two adjacent vertices have the same color.
A necessary and sufficient condition for a planar map to be bipartite is that all its faces have even degrees. 
This paper deals exclusively with planar bipartite maps.

Ref.~\IRRED\ was mainly devoted to the enumeration of {\it maps with a single boundary}, i.e.\ maps with a distinguished oriented 
edge (the root edge) and with a control on its {\it outer degree}, i.e.\ the degree of the face lying on the right of this edge (the root face).
By first choosing the root face, then the root edge, maps with a single boundary may alternatively be defined as {\it maps 
with a distinguished face of controlled degree} and with a marked oriented edge incident to this face
whose orientation is such that the distinguished face lies on its right. In this paper, we shall extend this definition to
maps with a number $r\geq 2$ of boundaries, namely maps with $r$ distinguished faces of controlled degrees, 
and with a marked oriented edge incident to each distinguished face (again with an orientation such that the distinguished face at hand
lies on its right). The faces which have not been distinguished will be referred to as {\it inner faces}. Note that
maps with $r\geq 2$ boundaries have at least two faces, so they cannot be reduced to trees. The case of
maps with two boundaries was actually discussed already at the end of Ref.~\IRRED\ and we shall rely here on the
corresponding results. 

Throughout this paper, we consider $d$ some even integer larger than or equal to 2.
We set for convenience
\eqn\db{d=2b}
with $b$ an integer larger than or equal to 1.
The results presented here concern the enumeration of so-called {\it $d$-irreducible} maps, as defined
now. Recall that the {\it girth} of a map is the minimal length (number of edges) of its cycles (simple closed paths). 
Note that the girth of a bipartite map is necessarily even.
A map with one or several boundaries is said {\it $d$-irreducible} if its girth is at least
$d$ and if all its cycles of length $d$ are the boundary of an inner face of degree $d$. 
Clearly, from the girth condition, all the faces in a $d$-irreducible map with $r\geq 2$ boundaries
have a degree larger than or equal to $d$. In this paper, we shall restrict our enumeration to the case where all 
the marked faces have a degree {\it strictly} larger than $d$ (i.e\ only inner faces may have degree $d$).
We shall enumerate $d$-irreducible maps with a weight $z$ per $d$-valent inner face and a 
weight $x_{2j}$ per $2j$-valent inner face for $j>b$ (recall that all faces have even degree in a bipartite map),
the marked faces receiving no weight.
We then denote by 
\eqn\Fdef{F_{2j_1,2j_2,\ldots,2j_r}^{(d)}(z;x_{d+2},x_{d+4},\ldots )}
 the generating function of $d$-irreducible maps with 
$r$ boundaries whose marked faces have respective degrees $2j_1, 2j_2,\ldots,2j_r$, with $j_\ell >b$, $\ell=1,\ldots r$.
The main result of this paper is a general expression for this generating function (see eq.~(2.14) below). It extends a formula due to Collet and Fusy \CF\ for the case of arbitrary (non necessarily irreducible) bipartite planar maps, which we recover by taking $d=0$ (see the discussion in section 4 below).
We shall then use this expression to derive a number of explicit formulas for the {\it numbers}
of  $d$-irreducible planar bipartite maps with prescribed numbers of $2j$-valent faces for each $j\geq b$ (section 3 below).

\newsec{Generating functions for $d$-irreducible maps}

\subsec{The starting point: $d$-irreducible maps with two boundaries}

Our starting point is a formula for the quantity \Fdef\ at $r=2$, namely the 
generating function $F^{(d)}_{2j_1,2j_2}$ of $d$-irreducible bipartite maps with 
two boundaries of respective lengths $2j_1$ and $2j_2$. From [\xref\IRRED, eq.~(9.17)], it is given by
 \eqn\Ftwobound{F^{(d)}_{2j_1,2j_2} = 2j_1 {2j_1-1 \choose j_1+b} 2j_2 {2j_2-1 \choose j_2+b}  {(R^{(d)})^{j_1+j_2} 
\over j_1+j_2} \qquad j_1,j_2 > b ,}
where the function $R^{(d)}$ is a fundamental object related via
\eqn\Rd{R^{(d)}=1+U^{(d)}_0}
to the first member $U^{(d)}_0$ of a family of {\it slice generating functions} $U_k^{(d)}$, $0\leq k\leq b$.
We shall not explain here what are precisely the ``slices" enumerated by $U_k^{(d)}$ and refer to \IRRED\ 
for such a discussion. We remind simply that slices appear as the result of a cutting of the
maps along geodesic paths, and have themselves a recursive decomposition.  As a consequence, 
the generating functions  $U_k^{(d)}$ are fully determined via a set of recursive equations
\eqn\recurU{U_k^{(d)}=z \delta_{k,b-1}+\sum_{q\geq 1} \sum_{
k_1,\ldots,k_q \geq 1 \atop k_1+\cdots+k_q=k+1 } \prod_{i=1}^{q} U_{k_i}^{(d)} \qquad 0\leq k\leq b-1 ,} 
which, together with \Rd, form a closed system, once completed by the relation 
\eqn\Ubd{U_b^{(d)} = \sum_{j \geq b+1} {2j-1 \choose j+b} x_{2j} (R^{(d)})^{j+b}.}
Our goal is to extend formula \Ftwobound\ to the case of maps with $r> 2$ boundaries.

\subsec{The tree representation}
\fig{Tree building rules in the case $b>1$ (see text).}{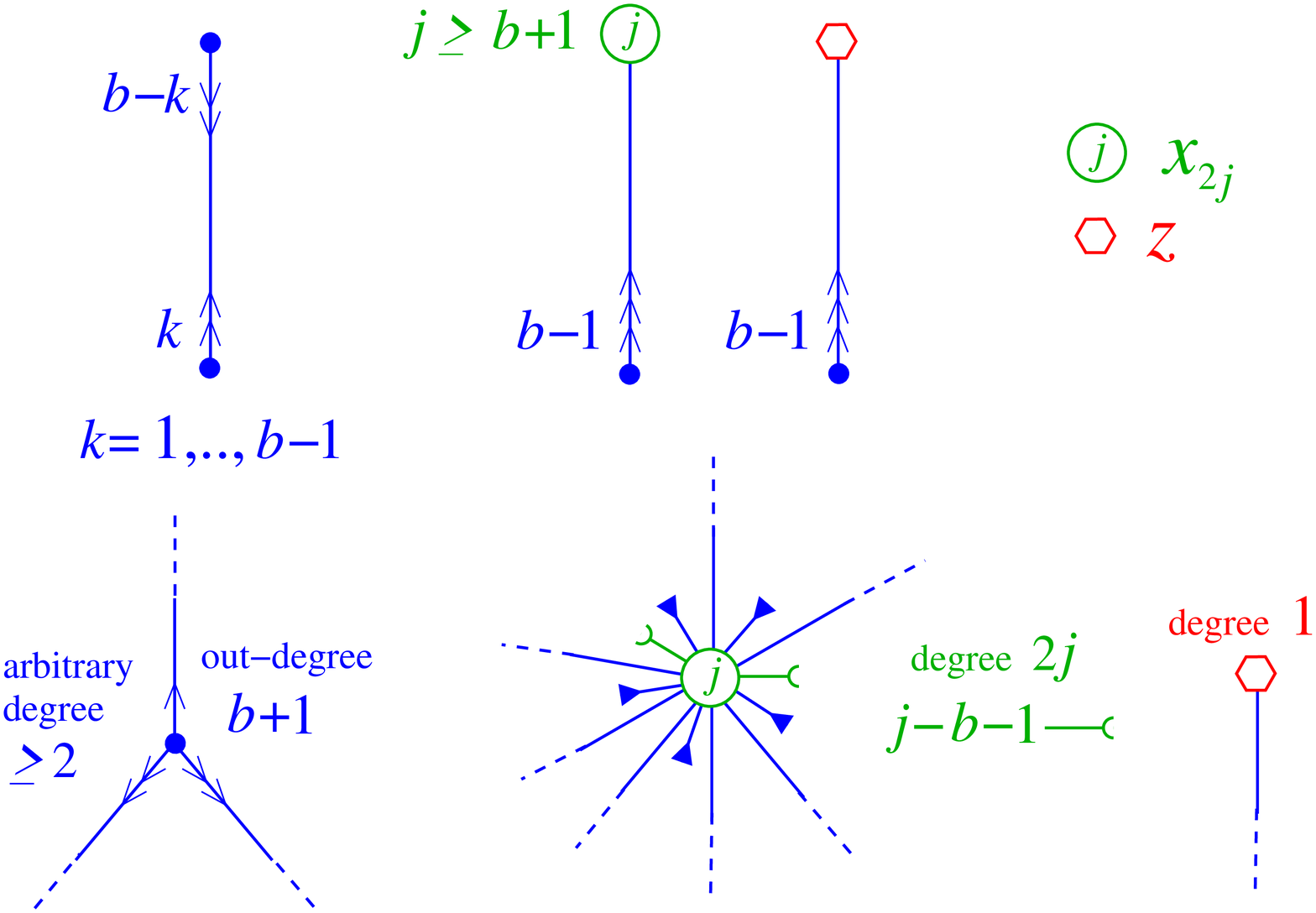}{11.cm}
\figlabel\trees
\fig{Interpretation of $R^{(d)}$, $U^{(d)}_k$, $1\leq k\leq b-1$ and $Z^{(d)}_j,$ $j\geq b+1$ as
generating functions for planted trees.
}{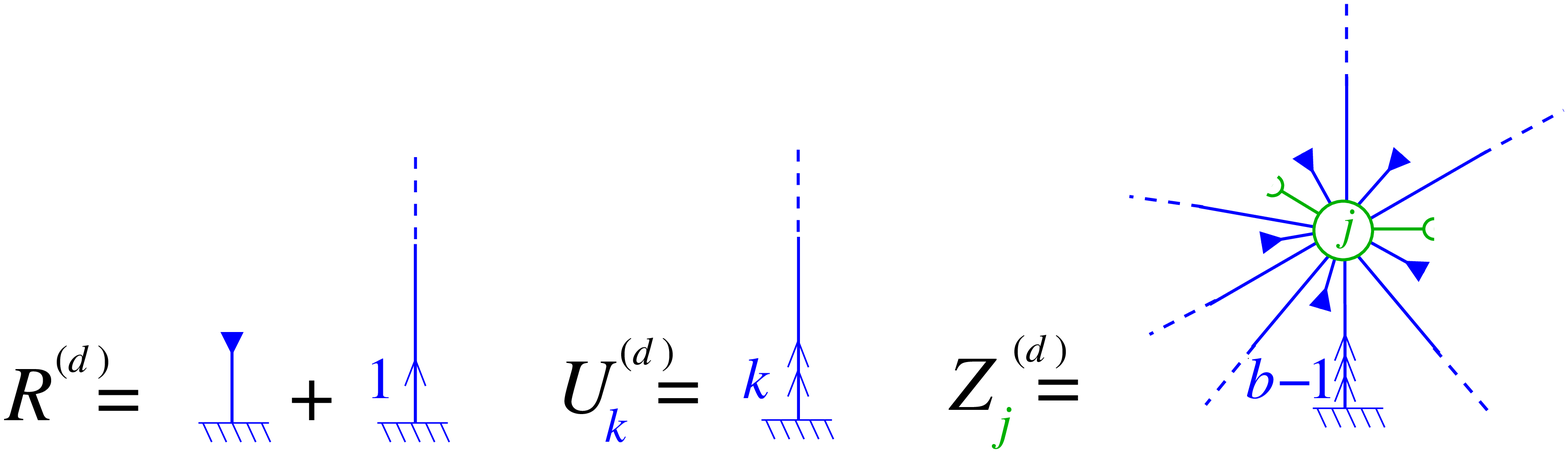}{12.cm}
\figlabel\gftrees
For $b>1$, we may eliminate explicitly $U^{(d)}_0$ and $U^{(d)}_b$ from the system \Rd-\Ubd\ by using $U^{(d)}_0=U^{(d)}_1$ 
(as obtained from \recurU\ at $k=0$) and using the expression \Ubd\ for  $U^{(d)}_b$, resulting in the system
\eqn\Rdrecu{\eqalign{R^{(d)}&=1+U_1^{(d)}\cr
U_k^{(d)}&=\sum_{q\geq 1} \sum_{
1\leq k_1,\ldots,k_q \leq b-1 \atop k_1+\cdots+k_q=k+1 } \prod_{i=1}^{q} U_{k_i}^{(d)} \qquad 1\leq k\leq b-2\cr
U_{b-1}^{(d)}&=\sum_{q\geq 2} \sum_{
1\leq k_1,\ldots,k_q \leq b-1 \atop k_1+\cdots+k_q=b } \prod_{i=1}^{q} U_{k_i}^{(d)}+\left(z + \sum_{j \geq b+1} Z^{(d)}_j\right)\cr
Z^{(d)}_j&= {2j-1 \choose j+b} x_{2j} (R^{(d)})^{j+b}\qquad j\geq b+1 .\cr }}
In view of this system, we may reinterpret $R^{(d)}$, $U_k^{(d)}$, $1\leq k \leq b-1$ and $Z_j^{(d)}$, $j\geq b+1$  as generating functions  for particular planted plane trees built according to the following rules, displayed in fig.~\trees. The trees are made of
\item{-}{four types of vertices:}
\itemitem{-}{black inner vertices (represented by filled dots);}
\itemitem{-}{white inner vertices (represented by open circles) carrying a label $j\geq b+1$;}
\itemitem{-}{leaf-vertices (represented by hexagons);}
\itemitem{-}{a univalent root-vertex (represented by the ground symbol in fig.~\gftrees);}
\item{-}{two types of edges:}
\itemitem{-}{bi-oriented edges of type $k/(b-k)$ with $1\leq k\leq b-1$, i.e.\ edges carrying $k\geq 1$ outgoing arrows pointing away 
from one of their extremities, and $b-k\geq 1$ outgoing arrows pointing away from the other extremity.
These bi-oriented edges connect only black inner vertices or the root;}
\itemitem{-}{mono-oriented edges carrying $b-1$ outgoing arrows pointing away 
from one of their extremities, being necessarily a black inner vertex or the root. The other
extremity (without arrows) is necessarily a white inner vertex or a leaf-vertex;}
\item{-}{with the vertex constraints:}
\itemitem{-}{black inner vertices have arbitrary degrees larger than or equal to $2$ but a fixed {\it out-degree}
equal to $b+1$. By out-degree of a vertex, we mean the total number of outgoing arrows pointing away from it. (Note
that the out-degree constraint restricts in practice the degree of black inner vertices, which can be at most $b+1$);}
\itemitem{-}{white inner vertices are incident to mono-oriented edges and to two types of decorations:
buds (represented with a semi-circular endpoint) and blossoms (represented with a triangular endpoint).
A white inner vertex with label $j$ has total degree $2j$ and a number $j-b-1$ of buds.
It is weighted by $x_{2j}$;}
\itemitem{-} leaf-vertices have degree $1$ and are weighted by $z$.
\par
\noindent The reader will easily check that the rules above are designed so as to reproduce the
equations \Rdrecu\ by a {\it canonical decomposition} of the (supposedly planted) trees into descending subtrees 
by cutting them at the level of their first vertex. This vertex can be a black inner vertex of degree $q+1$ ($q\geq 1$),
reproducing the $q$-dependent terms in \Rdrecu, a white inner vertex with label $j$, reproducing the 
$j$-dependent terms in \Rdrecu, or a leaf-vertex, reproducing the $z$-dependent term in \Rdrecu.
More precisely, $R^{(d)}$, $U_k^{(d)}$, $1\leq k \leq b-1$ and $Z_j^{(d)}$, $j\geq b+1$ are easily 
identified with the generating functions for the planted trees displayed in fig.~\gftrees.
Let us mention that similar trees (without leaf-vertices) were introduced in [\xref\BFa,\xref\BFb] in the context of maps with
controlled girth. This is consistent with the fact that, as observed in \IRRED, enumerating such maps indeed corresponds 
to imposing $z=0$ in our setting. 
\fig{Tree building rules in the case $b=1$.}{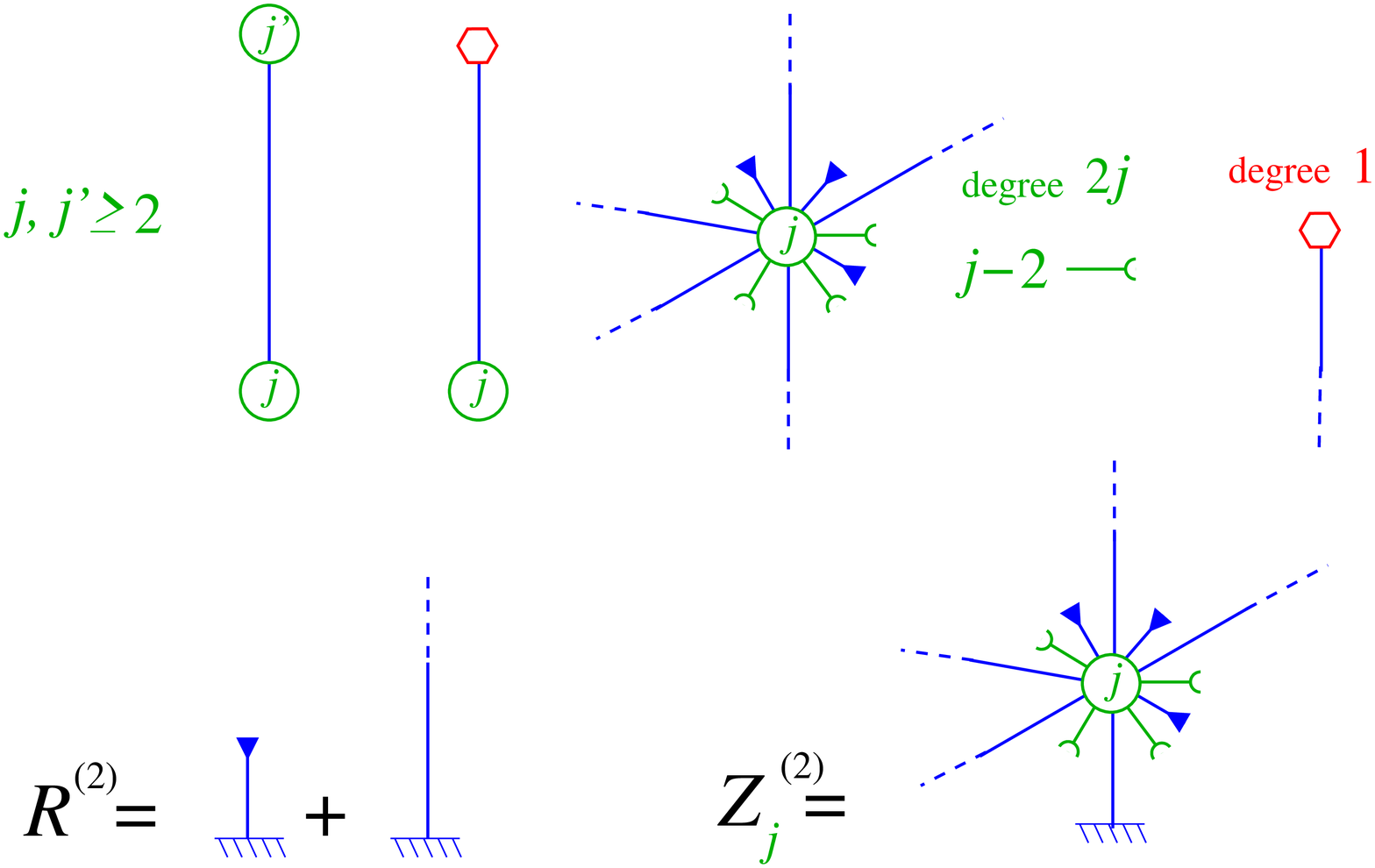}{12.cm}
\figlabel\treesbone
For $b=1$, eqs.~\Rd-\Ubd\  simply reduce to
\eqn\Rbone{\eqalign{R^{(2)}&=1+z+\sum_{j \geq 2} Z^{(2)}_j\cr
Z^{(2)}_j&= {2j-1 \choose j+1} x_{2j} (R^{(2)})^{j+1} \qquad j\geq 2\cr }}
and may be viewed as generating the simpler trees of fig.~\treesbone\ with un-oriented edges.

\subsec{A formula for the generating function of $d$-irreducible maps with $r\geq 2$ boundaries}

Returning to the general case $b\geq 1$, we now consider the generating function $ F^{(d)}_{2j_1,2j_2,\ldots,2j_r}$
for maps with  $r\geq 2$ boundaries. Since $r\geq 2$, maps enumerated by $ F^{(d)}_{2j_1,2j_2,\ldots,2j_r}$ may be obtained from 
maps enumerated by $F^{(d)}_{2j_1,2j_2}$ by  a simple marking of $r-2$ extra faces of respective degrees $2j_3,2j_4,\ldots 2j_r$.
After marking, these extra faces are no longer considered as inner faces. Still the constraint that $d$-cycles must be boundaries 
of inner faces of degree $d$ is not affected by the marking since the marked faces are all assumed to have degree strictly larger than $d$
(i.e.\ we assume $j_\ell>d$, $\ell=1,\ldots,r$).
At the level of generating function, the marking is performed via the action of derivatives with respect to $x_{2j_\ell}$ for $\ell=3,\ldots r$,
namely
\eqn\multibound{F^{(d)}_{2j_1,2j_2,\ldots,2j_r}=\left(\prod _{\ell=3}^r 2 j_\ell {\partial\over \partial x_{2j_\ell}}\right) F^{(d)}_{2j_1,2j_2}\ ,
\qquad r\geq 2}
(with the usual convention that the empty product represents the identity operator). 
Here the $2j_\ell$ factors account for the choice of an oriented edge incident to each newly marked face.
\fig{The building rule of the special vertex labeled by $j_1$ and $j_2$. It has degree $j_1+j_2+1$ 
and is incident to mono-oriented (or un-oriented if $b=1$) edges or blossoms only. The edge connecting it to the root is a 
mono-oriented edge (with thus $b-1$ arrows pointing away from the root -- it is un-oriented if $b=1$).}{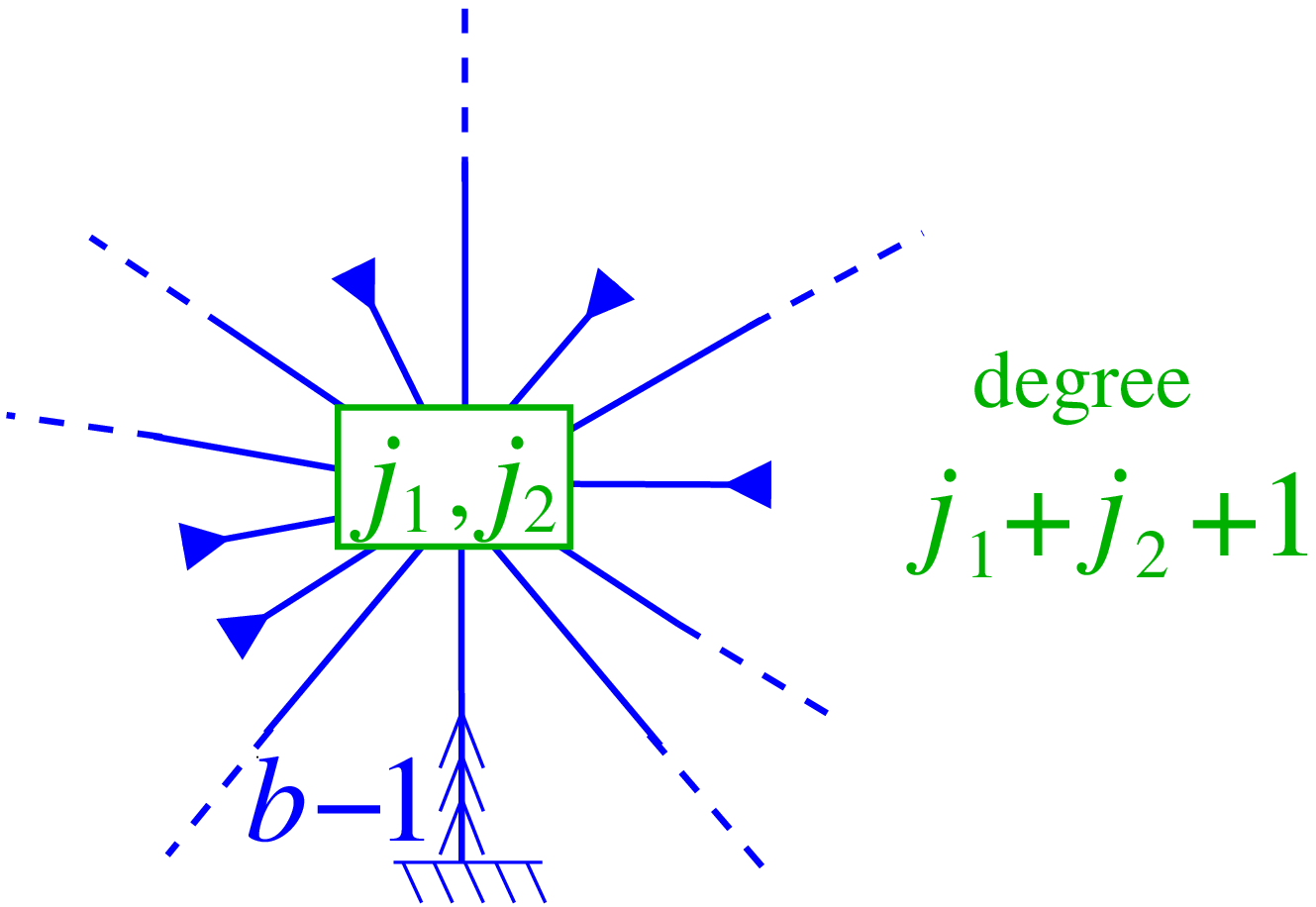}{7.cm}
\figlabel\specialroot
From \Ftwobound, we may thus write
\eqn\combin{F^{(d)}_{2j_1,2j_2,\ldots,2j_r}=\prod_{\ell=1}^r 2 j_\ell \ \times {{2j_1-1 \choose j_1+b} {2j_2-1\choose j_2+b}\over j_1+j_2}
\times \ H^{(d)}_{2j_1,2j_2;2j_3,\ldots,2j_r}}
where
\eqn\Gdef{H^{(d)}_{2j_1,2j_2;2j_3,\ldots,2j_r}= \left(\prod_{\ell=3}^r {\partial\over \partial x_{2j_\ell}} \right) 
(R^{(d)})^{j_1+j_2}\ .
}
In the tree language, the function $H^{(d)}_{2j_1,2j_2;2j_3,\ldots,2j_r}$ can be interpreted as the generating function for planted
trees built according to the rules of fig.~\trees\ (or of fig.~\treesbone\ if $b=1$), whose root is adjacent to a first special vertex carrying
both labels $j_1$ and $j_2$, as 
displayed in fig.~\specialroot, and with $r-2$ marked white inner vertices of respective degrees $2j_3, 2j_3,\ldots,2 j_r$.
The building rule of the special vertex is designed to reproduce the term $(R^{(d)})^{j_1+j_2}$ in \Gdef: this special 
vertex must have total degree $j_1+j_2+1$ 
and be incident to mono-oriented (or un-oriented if $b=1$) edges or blossoms only. The special vertex and the $r-2$ 
marked white inner vertices receive no weight.
Let us now show that we may write the generating function $H^{(d)}_{2j_1,2j_2;2j_3,\ldots,2j_r}$ as 
\eqn\resone{\eqalign{H^{(d)}_{2j_1,2j_2;2j_3,\ldots,2j_r}& =\prod_{\ell=3}^r {2j_\ell-1\choose j_\ell+b} \cr & \times 
\sum_{p=1}^{r-2} {r-3 \choose p-1} {\partial^p \over \partial z^p} (R^{(d)})^{j_1+j_2}  {\partial^{r-p-2} \over \partial z^{r-p-2}} (R^{(d)})^{\sum\limits_{\ell=3}^r (j_\ell+b)}\ . \cr}}
\fig{Schematic picture of the bijection between the trees enumerated by $H^{(d)}_{2j_1,2j_2;2j_3,\ldots,2j_r}$ with
a number $p$ of first generation vertices, supposedly numbered $1$ to $p$ is all possible ways (hence
the factor $p!$) and pairs made of (i) a forest of $p$ trees with roots labeled $1$ to $p$, with $(r-p-2)$ marked non-root white vertices 
and (ii) a tree with first special vertex and $p$ marked leaf-vertices labeled $1$ to $p$.
Here $p=3$ and $r=8$.}{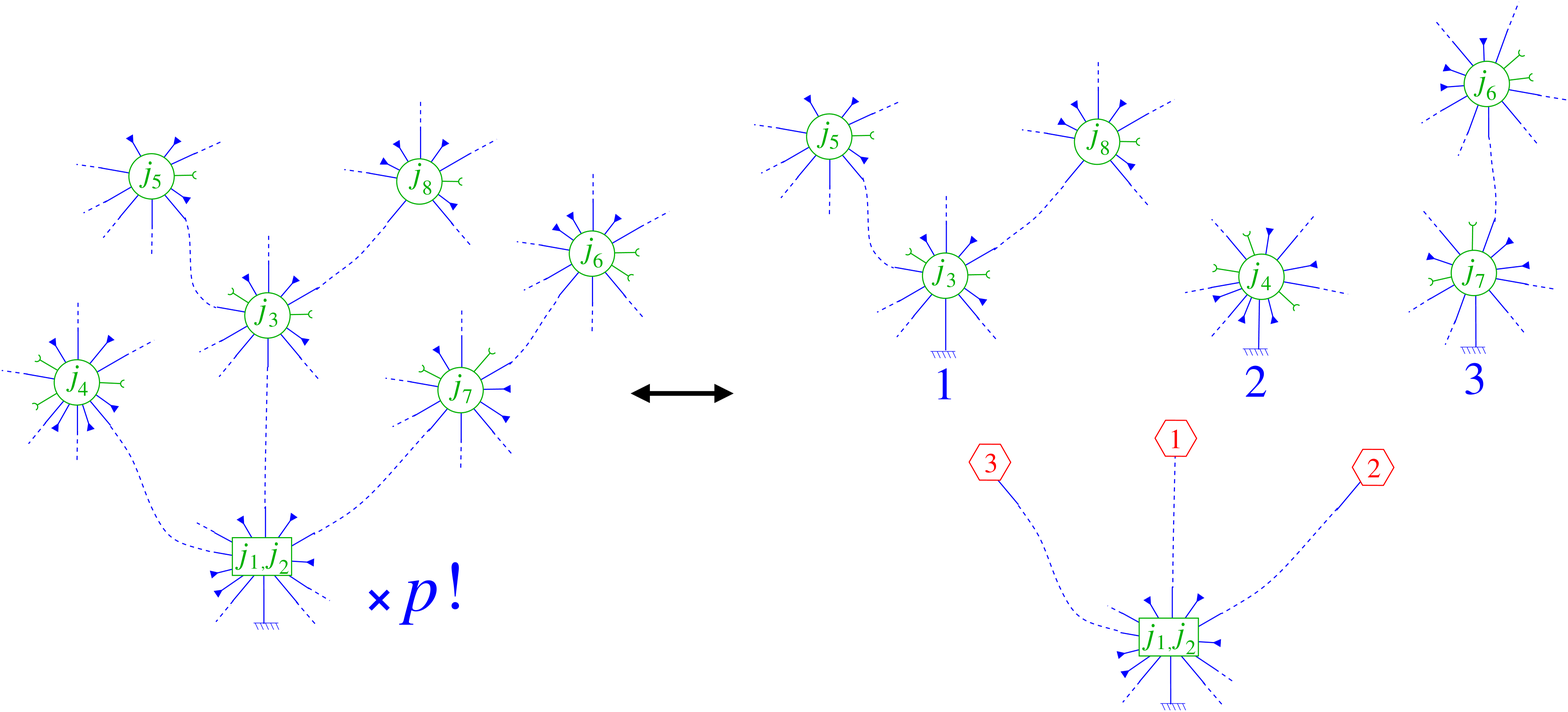}{14.cm}
\figlabel\foresttwo
We shall proceed in two steps: the trees enumerated by $H^{(d)}_{2j_1,2j_2;2j_3,\ldots,2j_r}$ have a number
$p\geq 1$ of ``first generation" marked white vertices, which have the special vertex as ``direct" ancestor, i.e.\ are such
that the branch from the special vertex to them does not pass via any other marked 
white vertex. We number these first generation vertices from $1$ to $p$ (in all possible ways so that each tree is counted $p!$ times -
see fig.~\foresttwo, left) and cut the tree at the level of their mono-oriented incident edge in the ascendent part of the tree. The $p$ cut parts form 
a forest made of $p$ rooted trees, whose roots are incident to the $p$ first generation vertices, and are labeled $1$ to $p$
(see fig.~\foresttwo, right).  This forest is equipped with a total of $(r-p-2)$ marked white vertices (those which are not
first generation). 
As for the part containing the special vertex, we may repair it by adding a leaf-vertex to each of the
cut mono-oriented edges, thus re-creating a tree satisfying the rules of fig.~\trees, with root adjacent to the special vertex, and with 
now $p$ marked (and labeled) leaf-vertices. These marked leaf-vertices receive no weight, so
the generating function for such trees is simply
\eqn\firstblock{{\partial^p\over \partial z^p} (R^{(d)})^{j_1+j_2}\ .}
This explains the $j_1,j_2$-dependent term in \resone. Moreover, the $p!$ to $1$ correspondence
displayed in fig.~\foresttwo\ can be made $1$ to $1$ by simply removing the labeling of the 
tree roots in the forest part.
\fig{Schematic picture of the correspondence between forests of $p$ unlabeled trees with $(r-p-2)$ marked
non-root white vertices and forests with $(r-2)$ unlabeled trees equipped with a total of $(r-p-2)$ marked leaf-vertices
labeled $1$ to $(r-p-2)$. The correspondence is $(r-p-2)!$ to $(r-3)!/(p-1)!$ (see text). Here $p=3$ and $r=8$.}{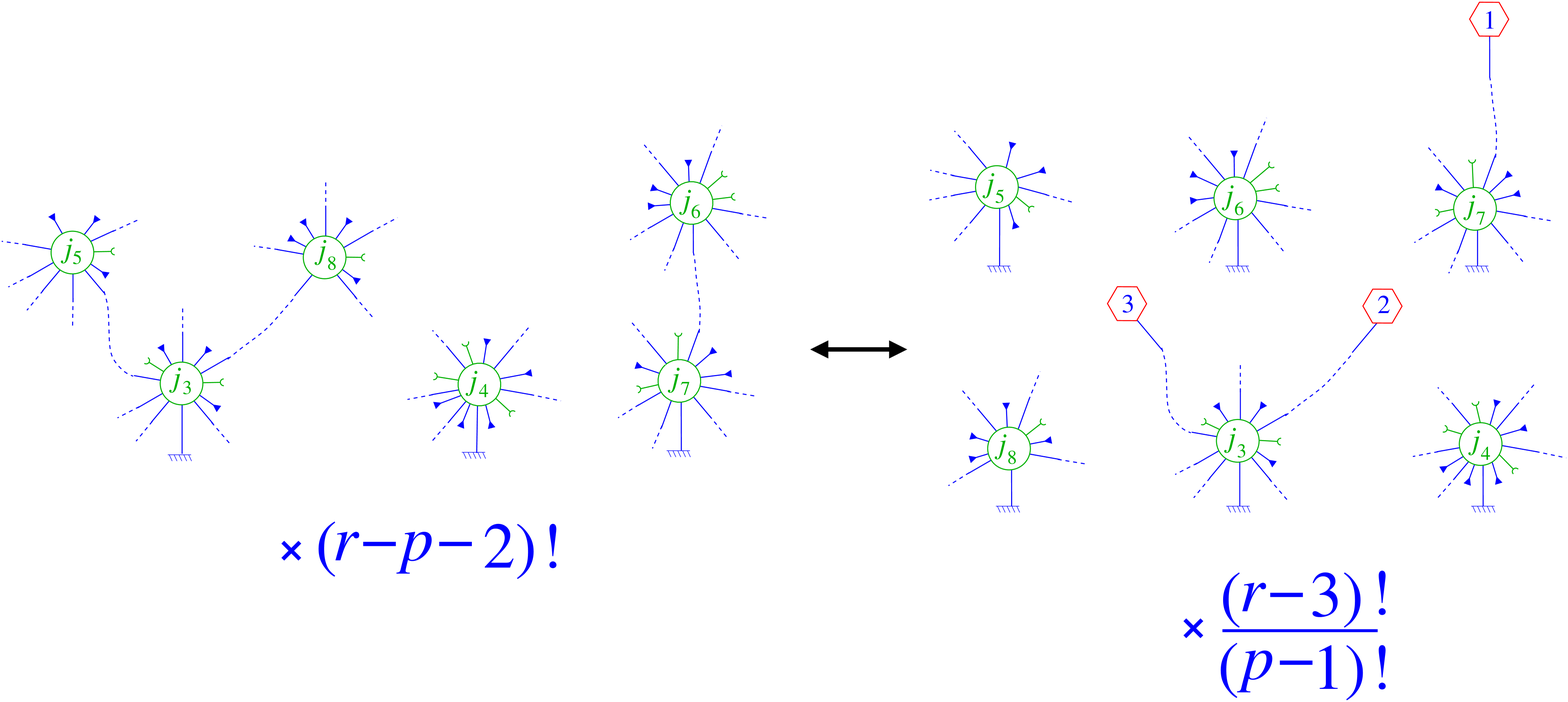}{14.cm}
\figlabel\forestone
In a second step, the enumeration of the complementary forest part is achieved as in \CF\ via a (many-to-many) correspondence between
forests of $p$ trees with a total of $(r-p-2)$ marked white non-root vertices (left side of fig.~\forestone) and 
forests of $(r-2)$ trees and a total of
$(r-p-2)$ marked leaf-vertices (right side of fig.~\forestone). We proceed as before by first 
labeling the $(r-p-2)$ marked white non-root vertices in all possible ways (each forest is thus counted $(r-p-2)!$ times) and cutting
the trees in the forest at the level of the mono-oriented edges leading to all these vertices. The cut mono-oriented
edges are as before completed by new added leaf-vertices in the {\it ascending part} of the trees and serve as roots
for the descending parts. The net result is a forest of $(r-2)$ trees whose roots are incident to the $(r-2)$ marked white vertices
with degrees $2j_3,\ldots,2j_r$, equipped with a total of $(r-p-2)$ marked leaf-vertices labeled $1$ to $(r-p-2)$. 
These forests are counted by 
\eqn\secondblock{{\partial^{r-p-2} \over \partial z^{r-p-2}}\prod_{\ell=3}^r \left( {2j_\ell-1\choose j_\ell+b} 
(R^{(d)})^{j_\ell+b} \right)} 
which explains the $(j_3\to j_r)$-dependent terms in \resone. The correspondence is not simply $(r-p-2)!$ to $1$ since, in 
\secondblock, we erased the information on how to reassemble the trees. The appropriate counting is more transparent if 
we start conversely from a forest of $(r-2)$ trees with a total of $(r-p-2)$ marked leaves: we then rebuild a forest of $p$ trees
and $(r-p-2)$ marked (and labeled) white non-root vertices by first replacing the marked leaf-vertex labeled $1$ by 
a descending subtree chosen among any of the trees in the forest which do not carry this marked leaf-vertex ($r-3$ choices),
creating a new forest {\it with one less tree}, then replacing the marked leaf-vertex labeled $2$ by 
a descending subtree formed by any of the trees in the new forest which do not carry this marked leaf-vertex ($r-4$ choices),
and continue until the $(r-p-2)$-th marked leaf-vertex in replaced by a descending subtree ($p$ choices). This leaves us with the desired forest 
with $p$ trees. Clearly, each configuration is counted $(r-3)!/(p-1)!$ times so that the correspondence is eventually
$(r-p-2)!$ to $(r-3)!/(p-1)!$, which explains the binomial factor ${r-3 \choose p-1}$ in \resone.

Now the sum in \resone\ may be simply evaluated through
\eqn\sumeval{\eqalign{& 
\sum_{p=1}^{r-2} {r-3 \choose p-1} {\partial^p \over \partial z^p} (R^{(d)})^{j_1+j_2}  {\partial^{r-p-2} \over \partial z^{r-p-2}} (R^{(d)})^{\sum\limits_{\ell=3}^r (j_\ell+b)}\cr &\quad = {\partial \over \partial \zeta}\left\{
\sum_{p=1}^{r-2} {r-3 \choose p-1} {\partial^{p-1} \over \partial z^{p-1}} (R^{(d)}(z+\zeta))^{j_1+j_2}  {\partial^{r-p-2} \over \partial z^{r-p-2}} (R^{(d)}(z))^{\sum\limits_{\ell=3}^r (j_\ell+b)} \right\}_{\zeta=0}
\cr &\quad = {\partial \over \partial \zeta }\left\{ {\partial^{r-3}\over \partial z^{r-3}} \left(
(R^{(d)}(z+\zeta))^{j_1+j_2} (R^{(d)}(z))^{\sum\limits_{\ell=3}^r (j_\ell+b)} \right)\right\}_{\zeta=0}\cr
 &\quad = {\partial^{r-3}\over \partial z^{r-3}}
 {\partial \over \partial \zeta }\left\{  (R^{(d)}(z+\zeta))^{j_1+j_2} (R^{(d)}(z))^{\sum\limits_{\ell=3}^r (j_\ell+b)}\right\}_{\zeta=0}
 \cr 
  &\quad = {\partial^{r-3}\over \partial z^{r-3}} {j_1+j_2 \over j_1+j_2+\sum\limits_{\ell=3}^r (j_\ell+b)}
 {\partial \over \partial \zeta }\left\{ (R^{(d)}(z+\zeta))^{{j_1+j_2}+\sum\limits_{\ell=3}^r (j_\ell+b)}\right\}_{\zeta=0}
 \cr 
 &\quad =  {j_1+j_2\over (r-2)b+\sum\limits_{\ell=1}^r j_\ell} {\partial^{r-2}\over \partial z^{r-2}}
 (R^{(d)})^{(r-2)b+\sum\limits_{\ell=1}^r j_\ell}\ .\cr
 }}
Here we simply used Leibniz formula to go from the second to the third line, as well as elementary operations.
Plugging this result in \resone\ and \combin,  we arrive at our main formula
\eqn\restwo{F^{(d)}_{2j_1,2j_2,\ldots,2j_r} ={1\over (r-2)b +
\sum\limits_{\ell=1}^r j_\ell} \prod_{\ell=1}^r 2 j_\ell {2j_\ell-1\choose j_\ell+b}  {\partial^{r-2}\over \partial z^{r-2}} (R^{(d)})^{(r-2)b+\sum\limits_{\ell=1}^r j_\ell}}
valid for $r\geq 2$ and $j_\ell \geq b+1$ for all $\ell=1,\ldots,r$. In this formula, the function $R^{(d)}$ is
the solution of eqs.~\Rd-\Ubd. As shown in \IRRED, this system reduces after elimination to the single equation for $R^{(d)}$
\eqn\eqforRd{z + \sum_{\ell=0}^{b} (-1)^{b-\ell} {b+\ell \choose 2\ell} {\rm Cat}(\ell) (R^{(d)})^{b-\ell}
+ \sum_{j \geq b+1} {2j-1 \choose j+b} x_{2j} (R^{(d)})^{b+j}=0}
which is algebraic if we impose an upper bound on the degree of the faces
(i.e.\ $x_{2j}$ vanishes for $j$ large enough).

\newsec{Enumeration formulas}

\subsec{A general formula}

A direct corollary of \restwo\ is a formula for the {\it number} $N^{(d)}_m(\{q_j\}_{j\geq b})$ of rooted bipartite planar $d$-irreducible maps
with outer degree $2m$ ($m\geq b+1$) and with $q_j$ faces of degree $2j$, $j\geq b$ (distinct from the root face if $j=m$). 
It is obtained by marking in the map {\it all the faces with degree strictly larger than $b$}, leaving the $d$-valent faces only as inner faces.
In \restwo, the function $R^{(d)}$ must therefore be evaluated at $x_{2j}=0$ for all $j>b$. We call $r^{(d)}(z)$ this function, which is obtained by solving \Rdrecu\ (or \Rbone\ if $b=1$)
at $Z^{(d)}_j=0$ for all $j$. It is determined by the equation
\eqn\eqforrd{z + \sum_{\ell=0}^{b} (-1)^{b-\ell} {b+\ell \choose 2\ell} {\rm Cat}(\ell) (r^{(d)})^{b-\ell} = 0,}
as read off eq.~\eqforRd\ at $x_{2j}=0$, $j>b$.

Eq.~\restwo\ translates directly into:
\eqn\microcanonical{\eqalign{N^{(d)}_m(\{q_j\}_{j\geq b})& = (2m){2m-1 \choose m+b}
\prod_{j>b} {1\over q_j!} {2j-1\choose j+b}^{q_j}\cr &\quad \times {1\over E+b(F-2-2q_b)} {(F-2)!\over q_b!} [z^{F-2}] (r^{(d)}(z))^{ E+b(F-2-2q_b)}\cr}}
where $F$ and $E$ are the total numbers of faces and edges respectively.
Note the absence of the factors $2j$ in the product which, in \restwo, accounted for a choice of edge incident
to each marked face. Here only the root face receives such a factor. Note also the factorial factors since,
as opposed to \restwo, the marked faces are not distinguished. 
We finally used the identities
\eqn\EF{\eqalign{F&= 1+q_b+\sum_{j>b}q_j\cr
E& =m+b\, q_b+\sum_{j>b}j\, q_j\cr}}
which lead to the identifications $r=1+\sum\limits_{j>b}q_j=F-q_b$ and $(r-2)b+\sum\limits_{\ell=1}^r j_\ell=
(r-2)b+m+\sum\limits_{j>b}j q_j=E+b(F-2-2q_b)$. Expression \microcanonical\ is valid for 
$F\geq 2+q_b$, i.e with at least two faces of degree strictly larger than $b$.

\subsec{Bipartite maps without multiple edges}
A first case of interest corresponds to bipartite maps without multiple edges, i.e.\ bipartite maps
of girth at least $4$. These maps are simply obtained from $2$-irreducible maps ($b=1$) by forbidding 
bivalent faces, i.e.\ setting $q_1=0$.
Using $r^{(2)}(z)=1+z$ (as seen directly from \Rbone\ at $Z^{(2)}_j=0$ for all $j$), and
\eqn\caldrtow{(F-2)! [z^{F-2}](1+z)^{E+F-2}={(E+F-2)!\over E!}}
we obtain the number ${\cal N}^{(4)}_m(\{q_j\}_{j\geq 2})$ ($= N^{(2)}_m(\{q_j\}_{j\geq 1})$ at $q_1=0$) of rooted bipartite planar 
maps {\it without multiple edges}, with outer degree $2m$ ($m\geq 2$) and with $q_j$ faces of degree $2j$, 
$j\geq 2$ (distinct from the root face if $j=m$):
\eqn\nnomult{{\cal N}^{(4)}_m(\{q_j\}_{j\geq 2}) = (2m){2m-1 \choose m+1}
\prod_{j>1} {1\over q_j!} {2j-1\choose j+1}^{q_j} {(E+F-3)!\over E!}}
with $E=m+\sum\limits_{j>1}j\, q_j$ the total number of edges, and $F=1+\sum\limits_{j>1}q_j$ the total number
of faces. This expression agrees with Proposition 35 in \BFb, up to a trivial rerooting factor.

For instance, we may obtain the number rooted $2m$-angulations ($m\geq 2$) without multiple edges, with $F$ faces,
by simply setting $q_j=(F-1)\delta_{j,m}$, namely
\eqn\twomang{2m {2m-1\choose m+1} {1\over (F-1)!} {2m-1\choose m+1}^{F-1} {(E+F-3)!\over E!}=
2m {2m-1\choose m+1}^F{((m+1)F-3)!\over (F-1)!(mF)!}}
since $E=m F$ in this case. For $m=2$, this yields the well-known number of quadrangulations without multiple
edges [\xref\TutteCPM-\xref\DSQ] 
\eqn\quadwithnomu{4 {(3F-3)! \over (F-1)!(2F)!}= 2 {(3F-3)! \over F!(2F-1)!}}
in terms of the number $F$ of faces. Note that our derivation requires at least two faces ($F\geq 2$) but
the above formula happens to also give the correct value $2$ at $F=1$.

\subsec{$4$-irreducible maps}
Another case of interest concerns $4$-irreducible maps, which requires the knowledge of 
$r^{(4)}(z)$. From \eqforrd, we immediately obtain
\eqn\valrfour{r^{(4)}(z)=1+{1-\sqrt{1-4z}\over 2}\ .}
After some straightforward algebra, we find, for $p>0$ and $s \geq 0$
\eqn\dqrfourp{\eqalign{{s! \over p}[z^{s}] (r^{(4)})^p
& = {\delta_{s,0} \over p} + (s-1)!
\sum_{\ell=0}^{\min(p-1,s-1) }{p-1 \choose \ell} {2(s-1)-\ell \choose s-1}\cr
&=  {\delta_{s,0} \over p} +{(2(s-1))!\over(s-1)! }\  {}_2F_1(1-p,1-s,2(1-s);-1)
\cr}}
in terms of the hypergeometric function ${}_2F_1$ (here the second term in the r.h.s. must be understood as
$0$ if $s=0$). Using this expression at $p=E+2F-4-4q_2$ and $s=F-2$,
we arrive at the formula
\eqn\Nfourl{\eqalign{N^{(4)}_m(\{q_j\}_{j\geq 2})& = (2m){2m-1 \choose m+2}
\prod_{j>2} {1\over q_j!} {2j-1\choose j+2}^{q_j}\cr &\ \times\!
\left\{{\delta_{F,2} \over E}\!
+\!{(2(F\!-\!3))!\over q_2! (F\!-\!3)! }\  {}_2F_1(5\!-\!E\!-\!2F\!+\!4q_2,3\!-\!F,2(3\!-\!F);-1)\right\}
\ ,\cr}}
valid for $F\geq 2+q_2$, i.e.\ with at least two faces of degree strictly larger than $4$ (note that setting $F=2$ 
implies necessarily $q_2=0$). As before, the second term being understood as $0$ if $F=2$. 
\fig{$4$-irreducible maps made of two hexagons and $q_2$ squares, for $q_2=0,1$ and $2$. For each
map, we indicated its multiplicity corresponding to the number of inequivalent possible choices of a root edge (not drawn here) 
incident to the outer hexagonal face.}{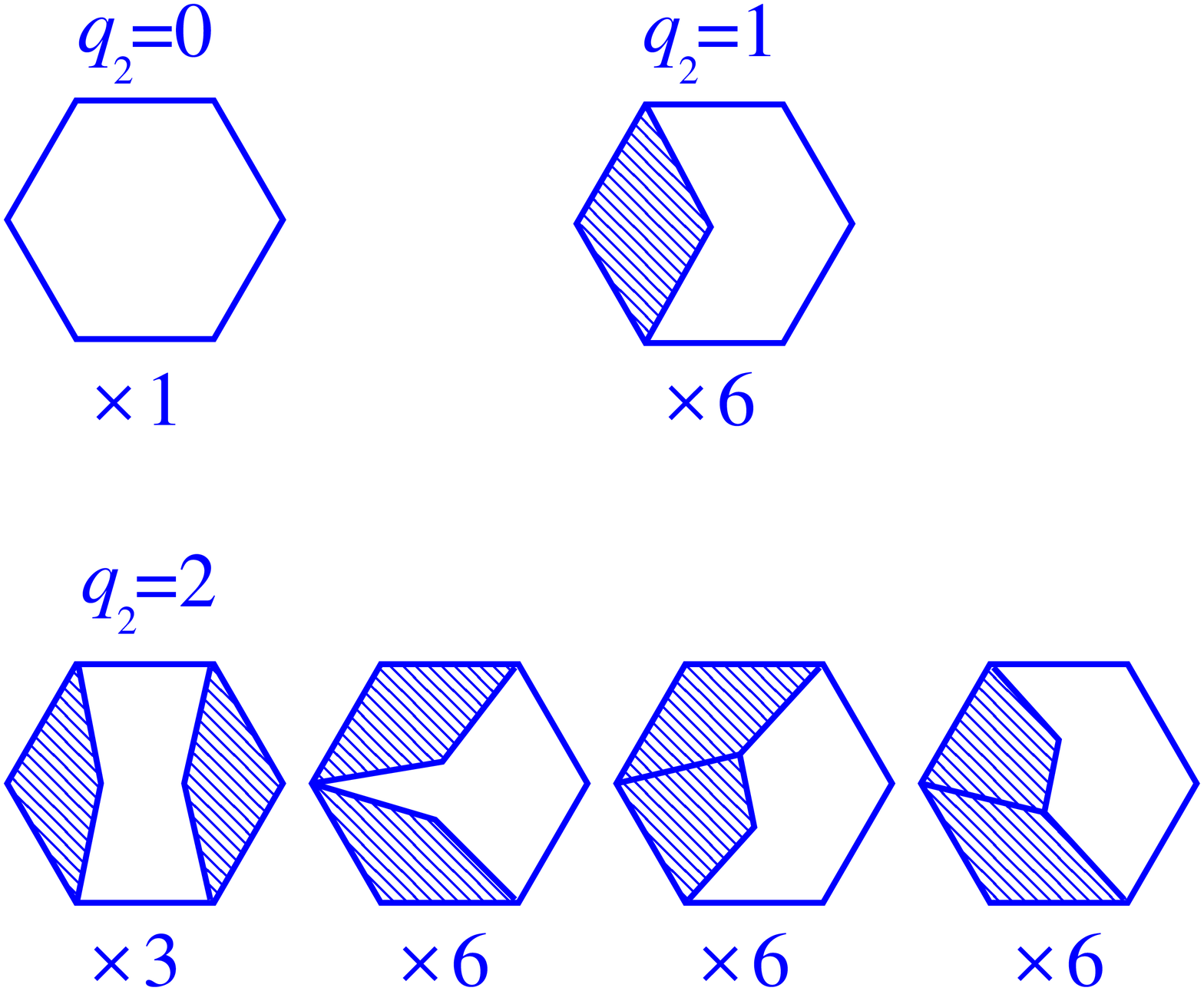}{8.cm}
\figlabel\sixsix
A simple application concerns maps with exactly two faces of degree strictly larger than $4$: the root face
of degree $2m$ ($m>2$) and another face of degree $2m'$ ($m'>2$). We have in this case
$q_j=\delta_{j,m'}+q_2 \delta_{j,2}$, $E=m+m'+2 q_2$ and $F=2+q_2$, so that we find a number
of rooted $4$-irreducible maps equal two
\eqn\fourmmprime{2m{2m-1\choose m+2} {2m'-1\choose m'+2} \left\{{\delta_{q_2,0}\over m\!+\!m'}
\!+\!{(2(q_2-1))!\over q_2!(q_2-1)!}\  {}_2F_1(1\!-\!(m\!+\!m'),1\!-\!q_2,2(1\!-\!q_2),-1) \right\}}
where the second term is to be understood as $0$ for $q_2=0$. For $m=m'=3$, it gives the 
sequence
\eqn\seq{1, 6, 21, 62, 180, 540, 1683, 5418, 17901, 60390, 207207, 720954, 2537964, \
9023328}
whose first three terms correspond to the maps displayed in fig.~\sixsix.
\subsec{Maps with girth at least $6$}
If we specialize the result of previous subsection to the case $q_2=0$, we obtain
the number of rooted maps with girth at least $6$, with outer degree $2m$ ($m\geq 3$) and with 
$q_j$ non-root faces of degree $2j$ ($j\geq 3$):
\eqn\girthsixl{\eqalign{{\cal N}^{(6)}_m(\{q_j\}_{j\geq 3})& = (2m){2m-1 \choose m+2}
\prod_{j>2} {1\over q_j!} {2j-1\choose j+2}^{q_j}\cr &\ \times\!
\left\{{\delta_{F,2} \over E}\!
+\!{(2(F\!-\!3))!\over (F\!-\!3)! }\  {}_2F_1(5\!-\!E\!-\!2F,3\!-\!F,2(3\!-\!F);-1)\right\}
\ . \cr}}
If all faces have degree $6$, i.e.\ $m=3$, $q_j=(F-1)\delta_{j,3}$, and $E=3F$, we obtain the
number of rooted hexangulations of girth $6$ with $F$ faces
\eqn\hexgirthsixl{
 \delta_{F,2} 
+6 \!{(2(F\!-\!3))!\over (F\!-\!1)!(F\!-\!3)! }\  {}_2F_1(5(1-\!F),3\!-\!F,2(3\!-\!F);-1)
}
(as before, the second term is to be understood as $0$ when $F=2$). It gives the sequence
\eqn\firsthex{
1, 3, 17, 128, 1131, 11070, 116317, 1287480, 14829188, 176250143, 2148687567}
corroborating the result of [\xref\BFa,\xref\OEIS].
 
\newsec{The case $b=0$}
It is interesting to include in our framework the case of general bipartite maps, without
constraints of irreducibility. We may indeed recover well-known formulas for bipartite maps by
simply extending our formulas to the case $b=0$. This should not come as a surprise since 
general bipartite maps correspond indeed to $0$-irreducible bipartite maps.  Eq.~\Ftwobound\ at
$b=0$ matches exactly the known formula \CF
 \eqn\Ftwoboundbone{F_{2j_1,2j_2} = 2j_1 {2j_1-1 \choose j_1} 2j_2 {2j_2-1 \choose j_2}  {R^{j_1+j_2} 
\over j_1+j_2} \qquad j_1,j_2 > b ,}
for the generating function of bipartite maps with two boundaries of length $2j_1$ and $2j_2$,
if we take as $R^{(0)}$ the solution $R$ of
\eqn\Rbzero{R =z+ \sum_{j \geq 1} Z_j, \qquad Z_j={2j-1 \choose j} x_{2j} R^{j}}
which replaces the system \Rd-\Ubd\ whenever $b=0$.
Here $z$ corresponds to a weight per {\it vertex} of the map, which again is not surprising since, in some sense, 
we may view the vertices as faces of degree $0$.
\fig{Building rules of the trees enumerated by $R$ and $Z_j$ satisfying eq.~\Rbzero.}{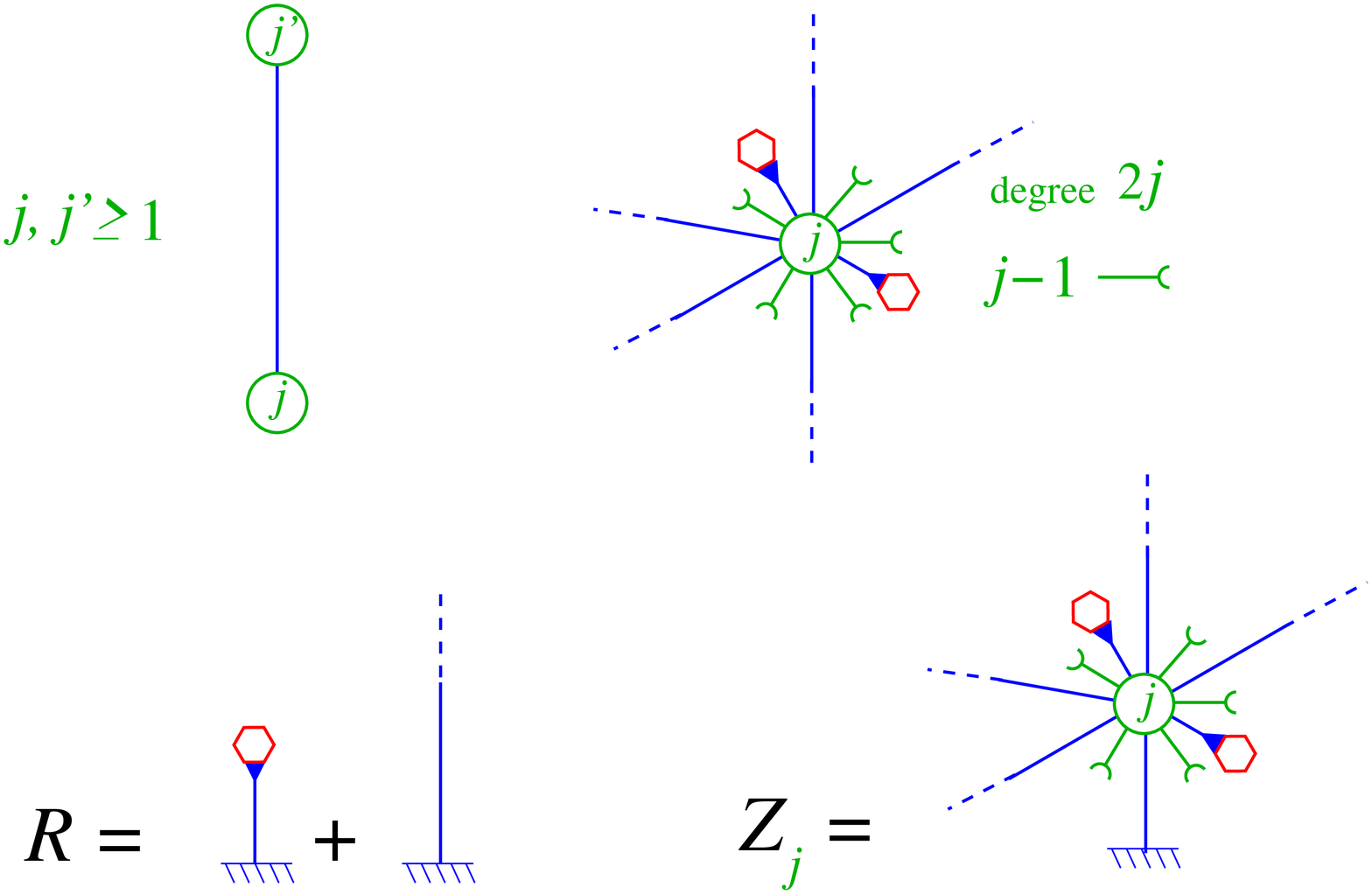}{12.cm}
\figlabel\treesbzero
The associated trees obey the rules displayed in fig.~\treesbzero, and are nothing but
the well-known ``blossom trees" introduced in [\xref\Schae,\xref\SchaeBCEM], with ``buds" and ``leaves" corresponding respectively 
to our buds and our blossoms enhanced by a leaf-vertex (see fig.~\treesbzero). 
Setting $b=0$ in \restwo\ yields the known beautiful formula \CF
\eqn\restwozero{F_{2j_1,2j_2,\ldots,2j_r} ={1\over
\sum\limits_{\ell=1}^r j_\ell} \prod_{\ell=1}^r 2 j_\ell {2j_\ell-1\choose j_\ell}  {\partial^{r-2}\over \partial z^{r-2}} R^{\sum\limits_{\ell=1}^r j_\ell}}
for the generating function of bipartite maps with multiple boundaries. 
Taking $r^{(0)}=z$ and formally $q_0=V$ and $F=2+E$ (indeed $F$ must be understood here as the number
of real faces plus that of degree $0$ faces -- i.e.\ vertices --, a sum which equals $2+E$ from Euler's relation), 
eq.~\microcanonical\ yields
\eqn\microcanonical{N_m(\{q_j\}_{j\geq 1}) = (2m){2m-1 \choose m}
\prod_{j\geq 1} {1\over q_j!} {2j-1\choose j}^{q_j}\ \times {(E-1)!\over V!} }
which is the well-known Tutte's formula for the number $N_m(\{q_j\}_{j\geq 1})$ of bipartite maps with root
face of degree $2m$ and with $q_j$ (non-root) faces of degree $2j$ \TutteCS. Here $E=m+\sum\limits_{j\geq 1} j q_j$
and $V=m+1+\sum\limits_{j\geq 1}(j-1)q_j$ denote the number of edges and vertices respectively .

\newsec{Concluding remarks}
Upon derivation with respect to $z$, formula \restwo\ yields 
\eqn\resthree{{\partial \ \over \partial z}F^{(d)}_{2j_1,2j_2,\ldots,2j_r} ={1\over (r-2)b+
\sum\limits_{\ell=1}^r j_\ell} \prod_{\ell=1}^r 2 j_\ell {2j_\ell-1\choose j_\ell+b}  {\partial^{r-1}\over \partial z^{r-1}} (R^{(d)})^{(r-2)b+\sum\limits_{\ell=1}^r j_\ell}\ .}
In this form, it remains valid at $r=1$ since, as shown in \IRRED, the generating function for $d$-irreducible maps
with a single boundary of length $2j_1$ satisfies the so-called ``pointing formula" 
\eqn\Fndpointing{{\partial F_{2j_1}^{(d)} \over \partial z} = {2j_1 \choose j_1-b} (R^{(d)})^{j_1-b}= {2j_1\over -b +j_1}{2j_1-1\choose j_1+b}
  (R^{(d)})^{-b+j_1} \qquad j_1> b .}
Conversely, we may use eq.~\Fndpointing\ as a starting point to derive eq.~\resthree\ from
the identity 
\eqn\multiboundone{F^{(d)}_{2j_1,2j_2,\ldots,2j_r}=\left(\prod _{\ell=2}^r 2 j_\ell {\partial\over \partial x_{2j_\ell}}\right) F^{(d)}_{2j_1}\ ,
\qquad r\geq 1}
by following the same procedure as above via a decomposition of the appropriate trees with marked white vertices into
forests. Still going back from \resthree\ to \restwo\  is not a simple matter as it requires integrating over $z$ and gives rise 
to some a priori unknown integration constant (which is $z$-independent but depends on the $j_\ell$'s). Fixing this constant 
requires knowing the value of $F^{(d)}_{2j_1,2j_2,\ldots,2j_r}$ at some particular value of $z$ (for instance at $z=0$ where it
enumerates maps with girth at least $d+2$), a problem whose difficulty might be comparable to that of the initial problem.
This is why we chose instead the generating function of maps with two boundaries as starting point.
 
Finally, we would like to stress that we used here trees as a simple tool to evaluate the wanted map generating functions
and did not recourse to any direct bijection between maps and trees. Nevertheless, we know from \IRRED\ that such direct bijections
do exist between the trees enumerated by $R^{(d)}$ and $U_k^{(d)}$, $k\geq 0$ and the so called slices or $k$-slices 
(which are particular instances of $d$-irreducible maps) enumerated by the same functions. Ref.~\IRRED\ describes in details a direct
correspondence between trees and slices in the case $x_{2j}=0$, $j>b$, and this correspondence can easily be extended to the case
of non-vanishing $x_{2j}$'s. To obtain a direct bijection between $d$-irreducible maps with several boundaries and trees,
we would simply need a direct bijection between $d$-irreducible maps with {\it two} boundaries of lengths $2j_1$ and $2j_2$
and the trees enumerated by the r.h.s of \Ftwobound. We have not found any such bijection so far. In particular, 
eq.~\Ftwobound\ still awaits a direct bijective proof. This is to be contrasted with eq.~\Fndpointing, which was given in \IRRED\
a direct bijective proof, which may be reformulated as a direct bijection with trees. However, as we just discussed, this formula is
not sufficient to recover our main result \restwo.

\bigskip
\noindent {\bf Acknowledgements:} The work of JB was partly supported by the ANR projects ``Cartaplus'' 12-JS02-001-01 and ``IComb'' ANR-08-JCJC-0011.
\listrefs
\end